\documentclass[12pt]{article}
\usepackage{amssymb}
\usepackage{amsmath}
\usepackage{fancybox}
\usepackage{tikz}
\usetikzlibrary{arrows,shapes,trees} % loads some tikz extensions

\newtheorem{theorem}{Theorem}[section]
\newtheorem{lemma}[theorem]{Lemma}

\newtheorem{definition}[theorem]{Definition}

\newcommand{\comments}[1]{}

\newcommand{\B}{\mathcal{B}}
\newcommand{\C}{\mathcal{C}}
\newcommand{\D}{\mathcal{D}}
\newcommand{\E}{{\rm E}}
\newcommand{\F}{\mathcal{F}}
\newcommand{\G}{\mathcal{G}}
\newcommand{\Ha}{\mathcal{H}}

\newcommand{\La}{{\rm La}}

\newcommand{\Oh}{\mathcal{O}}
\newcommand{\Pa}{\mathcal{P}}
\newcommand{\V}{\mathcal{V}}
\newcommand{\Sa}{\mathcal{S}}

\newcommand{\nchn}{\binom{n}{\lfloor \frac{n}{2}\rfloor}}
\newcommand{\mchm}{\binom{m}{\lfloor \frac{m}{2}\rfloor}}

\newcommand{\lanp}{\La(n,P)}

\newcommand{\ex}{{\rm Ex}}

\begin{document}
\title{On crown-free families of  subsets}
\author{
Linyuan Lu
\thanks{Department of Mathematics,
University of South Carolina, Columbia, SC 29208  USA ({\tt
lu@math.sc.edu}). This author was supported in part by NSF grant
DMS 1000475. } }
\maketitle

\begin{abstract}

  The crown $\Oh_{2t}$ is a height-2 poset whose Hasse diagram is a
  cycle of length $2t$.  A family $\F$ of subsets of
  $[n]:=\{1,2\ldots, n\}$ is {\em $\Oh_{2t}$-free} if $\Oh_{2t}$ is not
a weak subposet of  $(\F,\subseteq)$. Let $\La(n,\Oh_{2t})$ be the largest
  size of $\Oh_{2t}$-free families of subsets of $[n]$.  De~Bonis-Katona-Swanepoel proved $\La(n,\Oh_{4})=
{n\choose \lfloor \frac{n}{2} \rfloor} + 
{n\choose \lceil \frac{n}{2} \rceil}$.  Griggs and Lu proved that 
$\La(n,\Oh_{2t})=(1+o(1))\nchn$ for all even $t\ge 4$.  In this paper, 
 we prove  $\La(n,\Oh_{2t})=(1+o(1))\nchn$ for all odd $t\geq 7$.
\end{abstract}

\section{Introduction}
We are interested in estimating the maximum size of family of subsets
of the $n$-set $[n]:=\{1,\ldots,n\}$ avoiding a given (weak) subposet
$P$. The starting point of this kind of problem is Sperner's
Theorem from 1928~\cite{Spe}, which determined that the maximum size of an
antichain in the Boolean lattice $\B_n:=(2^{[n]},\subseteq)$ is $\nchn$.

For posets $P=(P,\le)$ and $P'=(P',\le')$, we say $P'$ is a {\em
weak subposet of $P$} if there exists an injection $f\colon P'\to P$
that preserves the partial ordering, meaning that whenever $u\le'
v$ in $P'$, we have $f(u)\le f(v)$ in $P$ (see \cite{Sta}). Throughout the
paper, when we say subposet, we mean weak subposet. The {\em
height\/} $h(P)$ of poset $P$ is the maximum size of any chain in
$P$.

A family $\F$ of subsets of $[n]$ can be viewed as a subposet of $\B_n$. If
$\F$ contains no subposet $P$, we say $\F$ is {\em $P$-free.}  We
are interested in determining the largest size of a $P$-free
family of subsets of $[n]$, denoted $\La(n,P)$.
 
In this notation, Sperner's Theorem~\cite{Spe} gives that
$\La(n,\Pa_2)=\nchn$, where $\Pa_k$ denotes the path poset on $k$
points, usually called a chain of size $k$.  Let $\B(n,k)$ be the
middle $k$ levels in the Boolean lattice $\B_n$ and
$\Sigma(n,k):=|\B(n,k)|$.  Erd\H{o}s \cite{Erdos45} proved that
$\La(n,\Pa_k)=\Sigma(n,k)$.  Griggs-Li-Lu \cite{GriLiLu} showed that
the similar results hold for a wide class of posets including diamonds
$\D_k$ ($A< B_1,\ldots, B_k <C$, for $k=3,4,7,8,9, 15,16,\ldots$),
harps $\Ha(l_1,l_2,\ldots,l_k)$ (consisting of chains
$\Pa_1,\ldots,\Pa_k$ with their top elements identified and their
bottom elements identified, for $l_1>l_2>\cdots>l_k$).

For any poset $P$, we define $e(P)$ to be the maximum $m$ such
that for all $n$, the union of the $m$ middle levels $\B(n,m)$
does not contain $P$ as a subposet. For any $\F\subset 2^{[n]}$,
define its Lubell value $h_n(\F):=\sum_{F\in \F}1/{n\choose |F|}$.
Let $\lambda_n(P)=\max\{h_n(\F)\colon \F\subset 2^{[n]}, P\mbox{-free}\}.$
A poset $P$ is  called {\em uniform-L-bounded} if $\lambda_n(P)\leq e(P)$
for all $n$. Griggs-Li \cite{GriLi2} proved $\La(n,P)=\Sigma(n, e(P))$ if
$P$ is uniform-L-bounded. The uniform-L-bounded posets include
$\Pa_k$ (for any $k\geq 1$),  diamonds $\D_k$ (for $k\in [2^{m-1}-1, 2^m-\mchm-1]$ where  $m:=\lceil \log_2(k+2)\rceil$), and harps  $\Ha(l_1,l_2,\ldots,l_k)$ (for $l_1>l_2>\cdots>l_k$), and other posets.

For any poset $P$, Griggs-Lu \cite{GriLu} conjectured the limit
$\pi(P):=\lim_{n\rightarrow\infty} \frac{\lanp}{\nchn}$ exists and
is  an integer. This conjecture is based on various known cases.
For example, an $r$-fork poset $\V_r$, which has elements $A<B_1,\ldots,B_r$, $r\ge2$.
Katona and Tarj\'an~\cite{KatTar} obtained bounds on $\La(n,\V_2)$ that
 he and DeBonis~\cite{DebKat}
 extended in  2007 to general $\V_r$, $r\ge2$, proving that
 \[
 \left(1+\frac{r-1}{n} + \Omega\left(\frac{1}{n^2}\right)\right)\nchn
 \le\La(n,\V_r)\le \left(1+ 2\frac{r-1}{n} +
 O\left(\frac{1}{n^2}\right)\right)\nchn.
 \]
 While the lower bound is
 strictly greater than $\nchn$, we see that $\La(n,\V_r)\sim \nchn$.
 Earlier, Thanh~\cite{Tha} had investigated the more general class
 of broom-like posets. Griggs and Lu~\cite{GriLu} studied the even
 more general class of baton posets.  These are tree posets
  (meaning that their Hasse diagrams are trees.)
 Griggs and Lu~\cite{GriLu} proved that $\pi(T)=1$ for any tree poset $T$
of height 2.  Bukh~\cite{Buk} proved that $\pi(T)=e(T)$
for any general tree poset $T$.

The most notable unsolved case is the diamond poset $D_2$.
Griggs and Lu first observed
$\pi(\D_2)\in [2,2.296]$.  Axenovich, Manske, and
Martin~\cite{AxeManMar} came up with a new approach which improves the
upper bound to $2.283$. Griggs, Li, and Lu \cite{GriLiLu} further
improves the upper bound to $2.27\dot 3=2\frac{3}{11}$. Very recently,
Kramer-Martin-Young \cite{KMY} recently proved $\pi(\D_2)\leq 2.25$.

The crown $\Oh_{2t}$ is another family of posets, which are neither
trees nor uniform-L-bounded. For $k\ge2$, the crown $\Oh_{2t}$
is a height-$2$ poset whose Hasse diagram is a cycle of length $2t$.
For $t=2$, $\Oh_4$ is also known as the butterfly poset; 
De Boinis-Katona-Swanepoel \cite{DebKatSwa} proved $\La(n,\Oh_{4})=\Sigma(n,2)$.
Griggs and Lu~\cite{GriLu} proved that $\La(n,\Oh_{2t})=(1+o(1))\nchn$ for all even $t\ge4$.
 For odd $t\ge 3$,  Griggs and Lu  showed that
$\La(n,\Oh_{2t})/\nchn$ is asymptotically at most
$1+\frac{1}{\sqrt2}$, which is less than 2. In this paper, we
determine all $\pi(\Oh_{2t})$ except for $\Oh_6$ and $\Oh_{10}$.

\begin{theorem}\label{crown}
 For odd $t\geq 7$, we have $\La(n, \Oh_{2t})=(1+o(1))\nchn$.
\end{theorem}
 
The proof of this theorem uses the concept of a {\em $k$-partite representation},
which was originally introduced by Conlon \cite{Conlon} to prove a similar
 Tur\'an-type result on hypercubes. (Conlon's result will be stated in Section 2.)
%We will review 
%Tur\'an problems on hypercubes in last section.
%Conlon's concept {\em $k$-partite representation} can be adapted for 
%posets of height $2$.

\begin{definition}\label{d1}
  A poset $P$ of height $2$ has a {\em $k$-partite representation} if
  there exist two integers $k$, $l$, and a family $\Pa\subseteq {[l]\choose
    k-1}\cup {[l]\choose k}$ such that
\begin{itemize}
\item  The poset $(\Pa, \subseteq)$ contains $P$ as a subposet.

\item And $G:=G(\Pa)$, a $k$-uniform hypergraph with $V(G)=[l]$ and $E(G)=\Pa\cap 
 {[l]\choose k}$ is $k$-partite.
%\item  there exists a fuction $\sigma\colon [l]\to [k]$ such that, 
% for each maximum element $\{i_1,\ldots, i_k\}$ of $\Pa$
% the image
% $\{\sigma(i_1),\ldots, \sigma(i_k)\}$ under $\sigma$ is $[k]$.
\end{itemize}
\end{definition}

%For example, $\Oh_{14}$ has a $3$-partite representation:
%with $l=7$, $k=3$, 
%$\Pa=$

%We actually proved the following analogue of Conlon's theorem.
Here is our main result.
\begin{theorem}\label{main}
  Suppose that a poset $P$ of height $2$ has a $k$-partitie representation for some $k\geq 2$.
Then $\La(n, P)=(1+o(1))\nchn$.
\end{theorem}

Conlon \cite{Conlon} proved that for all crowns $\Oh_{2t}$ except for $t=2,3,5$ 
have $k$-partitie representations for some $k$. For example,
$\Oh_{4t}$ (for $t\geq 2$) has a $2$-partitie representation $\Pa$
such that $G(\Pa)$ is the even-cycle $C_{2t}$. Similarly,
$\Oh_{2kt}$ (for $t\geq 2$) has a $k$-partitie representation $\Pa$
such that $G(\Pa)$ is the tight $k$-uniform cycle $C^{k}_{kt}$.
The first non-trivial case is  $\Oh_{14}$. 
The following $3$-representation of $\Oh_{14}$ is given by Conlon \cite{Conlon}: 
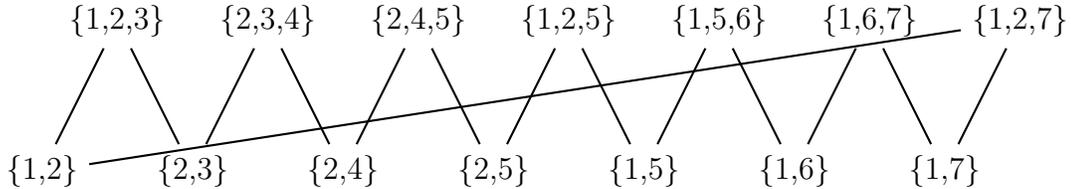
\begin{center}
\unitlength=1cm
  \begin{tikzpicture}[scale=1, thick]

% \foreach \place/\name/\string in 
%     {{(0,0)/1/$\{1,2\}$}, {(0.5,1)/2/$\{1,2,3\}$},
%      {(1,0)/3/$\{2,3\}$}, {(1.5,1)/4/$\{2,3,4\}$}}
%     {
%        \node (\name) at \place {};
%     }

\node (1) at (0,0) {\{1,2\}};
\node (2) at (1 ,2) {\{1,2,3\}};
\node (3) at (2,0) {\{2,3\}};
\node (4) at (3,2) {\{2,3,4\}};
\node (5) at (4,0) {\{2,4\}};
\node (6) at (5,2) {\{2,4,5\}};
\node (7) at (6,0) {\{2,5\}};
\node (8) at (7,2) {\{1,2,5\}};
\node (9) at  (8,0) {\{1,5\}};
\node (10) at (9,2) {\{1,5,6\}};
\node (11) at (10,0) {\{1,6\}};
\node (12) at (11,2) {\{1,6,7\}};
\node (13) at (12,0) {\{1,7\}};
\node (14) at (13,2) {\{1,2,7\}};

 \foreach \from/\to in {1/2, 2/3, 3/4, 4/5, 5/6, 6/7, 7/8,
   8/9, 9/10, 10/11, 11/12, 12/13, 13/14, 14/1}
 \draw  (\from) -- (\to);
  \end{tikzpicture}
\end{center}
Here $k=3$, $l=7$, and 
\begin{align*}
\Pa&=\{\{1,2\},  \{2,3\}, \{2,4\}, \{2,5\},\{1,5\},  \{1,6\},  \{1,7\},\\
&\{1,2,3\},  \{2,3,4\},  \{2,4,5\}, \{1,2,5\}, \{1,5,6\}, \{1,6,7\}, \{1,2,7\}\}.   
\end{align*}
It is easy to check that all the 3-edges $\{1,2,3\}$, $\{2,3,4\}$,
$\{2,4,5\}$, $\{1,2,5\}$, $\{1,5,6\}$, $\{1,6,7\}$, $\{1,2,7\}$ form a
$3$-partite $3$-uniform hypergraph. Thus, $\Pa$ is a $3$-partite representation of
$\Oh_{14}$.

For $t\geq 4$ and $t\not=5$, $\Oh_{2t}$ has a $k$-partite representation
for some $k$ (see \cite{Conlon}). It implies $\La(n, \Oh_{2t})=(1+o(1))\nchn$. 
Theorem \ref{crown} is a corollary of Theorem \ref{main}.
We also give an alternative proof for Griggs-Lu's result 
$\La(n, \Oh_{4t})=(1+o(1))\nchn$ for $t\geq 2$.

% Note $\Oh_{4}$ is the same as the butterfly poset $\B$; so $\La(n,\Oh_4)=\Sigma(n,2)
% =(2+o(1))\nchn$ by \cite{DebKatSwa}. Now the only unsolved Crowns are $\Oh_6$ and $\Oh_{10}$.

The rest of the paper is organized as follows. In section 2,
we will first review Conlon's theorem on Tur\'an problems on hypercubes; 
then we will prove an interesting Tu\'an-Ramsey result 
for $k$-partite $k$-uniform hypergraphs. Finally Theorem  \ref{main} 
will be proved in section 3.

% We will also list several lemmas from
%our preview paper \cite{GriLu}. We will prove Theorem \ref{crown}
%in section 3.

\section{Tru\'an problems on hypergraphs}
\subsection{ Tur\'an problem on hypercubes}
The problem of determining $\La(n, \Oh_{2t})$ is closely related to
the Tur\'an problem on the hypercube $Q_n$, i.e., the Hasse diagram of the Boolean
lattice $\B_n$. 
%Here the hypercube $Q_n$ is
%the graph on the vertex set of all binary strings of length $n$; a
%pair of two binary strings forms an edge of $Q_n$ if  they
%differs only on one bit. 
Erd\H{o}s \cite{Erdos64} first posed the problem of determine the size
of maximum subgraph of hypercube $Q_n$ forbidding a cycle
$C_{2k}$. Let $\ex(H, Q_n)$ be the maximum size of a subgraph of $Q_n$
forbidding a given graph $H$. Let $\pi(H, Q_n)=\lim_{n\to\infty}
\frac{\ex(H, Q_n)}{|E(Q_n)|}$.  This limit always exists.  Chung
\cite{chung} proved that $\pi(C_{4k}, Q_n)=0$ for all $k\geq 2$.  Alon
et al. \cite{AKT, ARSV} gave a characterization of all subgraphs $H$
of the hypercube which are Ramsey, that is, such that every
$k$-edge-colouring of a sufficiently large $Q_n$ contains a
monochromatic copy of $H$; in particular, $C_{4k+2}$ (for $k\geq 2$)
are Ramsey. F\"uredi and \"Ozkahya \cite{c14, c14+} showed that, for
$t > 3$, $\pi(C_{4t+2},Q_n)= 0$.  
Conlon \cite{Conlon} proved the following theorem, which covers all
known bipartite graphs $H$ with $\pi(H,Q_n)=0$.

\begin{theorem}[Conlon's Theorem \cite{Conlon}]
Suppose that $H$ is the Hasse diagram of a height-2 poset, which
admits a $k$-partite representation.  Then 
$\pi(H,Q_n)=0$. 
\end{theorem}
  In \cite{Conlon}, the $k$-partite
representation is defined over bipartite graphs. His
definition is equivalent to ours.  Conlon \cite{Conlon} observed
$C_{2t}$ (for $t\geq 4$ and $t\not=5$)
 admits a $k$-partite representation for some
$k$; thus, his result implies $\pi(C_{2t}, Q_n)=0$ for all 
$t\geq 4$ except for $t=5$.

\subsection{A Lemma on $k$-partitite $k$-uniform hypergraph}
Conlon \cite{Conlon}  used the following classical result of
Erd\H{o}s \cite{Erdos64} regarding the extremal number of complete
$k$-partite $k$-uniform hypergraphs.

\begin{lemma}\label{turan}
  Let $K^{(k)}_k(s_1,\ldots, s_k)$ be the complete $k$-partite
  $k$-uniform hypergraph with partite sets of size $s_1,\ldots,
  s_k$. Then any $K^{(k)}_k (s_1, \ldots, s_k)$-free $r$-uniform hypergraph
can have at most $O(n^{k-\delta})$ edges, where $\delta = \left(\prod^{k-1}_{i=1}
    s_i\right)^{-1}$.
\end{lemma}

In the scenario of the Boolean lattice, for any poset $P$ having
$k$-partite representation, we need prove that any family $\F$ of size
$(1+\epsilon)\nchn$ contains $P$.  Note that $\F$ is much sparser
comparing to the full Boolean lattice $2^{[n]}$. Lemma \ref{turan} is
not strong enough for our purpose.  We need the following lemma for
Ramsey-Tur\'an problems on hypergraphs, which may have independent
interest.

\begin{lemma}\label{cturan}
  For any positive integers $k$, $s_1,\ldots, s_k$, and $r$, consider a
  collection $\Ha:=\{H_i\}_{i\in I}$ (with an index set $I$)
of $k$-uniform hypergraphs over a
  common vertex set $[n]$.  Suppose that for each $i\in I$, $H_i$ does
  not contain $K^{(k)}_k (s_1, \ldots, s_k)$ as a sub-hypergraph, and
  for each $S\subset {[n]\choose k-1}$ there are at most $r$ hypergraphs $H_i$
 having
  edges containing $S$. Then, the total number of edges in this family
  is at most $ O(n^{k-\delta})$, where $\delta =
  \left(\prod^{k-1}_{i=1} s_i\right)^{-1}$.
\end{lemma}

{\bf Remark:} Since every hypergraph $H_i$ contains no $K^{(k)}_k
(s_1, \ldots, s_k)$, then $|E(H_i)|=O(n^{k-\delta})$ by Lemma
\ref{turan}.  This lemma says if the family of hypergraphs cover each
$(k-1)$-set at most $r$ times then the total number of edges is still
$O(n^{k-\delta})$, where the hidden constant in $O(\cdot)$ depends
on $k$, $s_1,\ldots, s_k$, and $r$, but not on $n$.

Our proof extensively uses the following convexity inequality, (also see
Lemma 2.3 of \cite{GriLu}.)
Suppose that $X$ is a random variable taking non-negative integer values.
If for any positive integer $s$, $\E(X)>s-1$, then
\begin{equation}
  \label{eq:convexity}
  \E {X\choose s}\geq {\E X \choose s}.
\end{equation}

% \begin{lemma}
% For any integer $s\geq 1$ and  any non-negative numbers $x_1,\ldots, x_m$, if 
% $\bar x:=\frac{\sum_{i=1}^m x_i}{m}\geq s-1$, then 
% $$\sum_{i=1}^m {x_i\choose s} \geq m {\bar x\choose s}.$$ 
% \end{lemma}
% {\bf Proof:} It is trival for $s=1$.
% Let $f(x):={x\choose s}$ if $x\geq s$ and $0$ otherwise.
% $f(x)$ is a well-defined convex function.

\noindent {\bf Proof of Lemma \ref{cturan}:} Let $H$ be the hypergraph
on the vertex set $[n]$ with $E(H)=\cup_{i\in I}E(H_i)$.  Observe that
each edge in $H$ can appear in at most $r$ $H_i$'s. Thus,
$$\sum_{i\in I}|E(H_i)|\leq r |E(H)|.$$
Since $r$ is a constant, it suffices to prove $|E(H)|=O(n^{k-\delta})$.
Deleting overlapped edges will not affect the magnitude of $|E(H)|$.
Without loss of generality, we can assume that edges of different $H_i$ are
distinct. If an edge $F$ of $H$ is in $H_i$, 
then we say this edge has color $i$. 
By hypothesis, $H$ has no monochromatic copies of $K^{(k)}_k(s_1,\ldots, s_k)$.

Without loss of generality, we assume $n$ is divisible by $k$ and write
$n=km$. Consider a random $k$-partition of $[n]=V_1\cup V_2\cup
\cdots \cup V_k$ where each part has the equal size $m$. We say an edge $F$ is 
{\em crossing} (to this partition), if
$F$ intersects every $V_i$ with exactly once. The probability of
an edge $F$ being crossing is
$$
\Pr(F \mbox{ is crossing})=\frac{\left(\frac{n}{k}\right)^k}{{n\choose k}}>
\frac{k!}{k^k}.$$
There exists a partition so that the number of crossing edges in $H$ 
at least $\frac{k!}{k^k}|E(H)|$.

Now we fix this partition $[n]=V_1\cup \cdots V_k$. Let $H'$ be the
subgraph consisting of all crossing edges in $H$ and
$H_i'$ be the subgraph consisting of all crossing edges in $H_i$ for $i\in I$.
It is sufficient to show $|E(H')|=O(m^{k-\delta})$, since $n=km$
and $k$ is a constant.

Set $|E(H')|\approx C m^{k-\delta}$ (with a big constant $C$ chosen later).
For $t_i\in \{1, s_i\}$ with $i=1,2\ldots, k$,
we would like to estimate the number of monochromatic (ordered)
copies, denoted by $f(t_1, t_2, \ldots, t_k)$, of
$K^{(k)}_k(t_1,\ldots, t_k)$ with the first $t_1$ vertices  in
$V_1$, the second $t_2$ vertices  in $V_2$, and so
on. 

{\bf Claim a:} For $0\leq l\leq k-1$, we have
$$f(s_1,\ldots, s_l, 1,\ldots, 1)\geq 
\left(1+o(1)\right)
\frac{\left(\frac{C}{m^\delta}\right)^{\prod_{j=1}^l s_j}} 
{\prod_{j=1}^l \left(s_j! r^{s_j-1}\right)^{\prod_{u=j+1}^ls_u}}
m^{k -l+ \sum_{j=1}^ls_j}.
$$ 

We prove claim (a) by induction on $l$.
For the initial case $l=0$, the claim is trivial since 
$f(1,1,\ldots, 1)=|E(H')|\approx C m^{k-\delta}$.

We assume Claim (a) holds for $l$. Now consider the case $l+1$. 
For any $S\in {V_1\choose s_1}\times \cdots \times {V_l\choose s_l} 
\times V_{l+2}\times \cdots \times V_k$, let $d_S^i$ be the 
number of vertices $v$ in $V_{l+1}$ such that all edges in the induced subgraph
of $H'$ on $S\times \{v\}$ have color $i$. Let $d_S=\sum_{i\in I}d_S^i$.
We have
\begin{align}
\label{eq:si}
f(s_1,\ldots, s_l,1,1, \ldots, 1) &= \sum_{S}\sum_{i\in I} d_S^i; \\
f(s_1,\ldots, s_l, s_{l+1},1,\ldots,1)&=\sum_{S}\sum_{i\in I} {d_S^i\choose s_{l+1}}.  
\label{eq:sil}
\end{align}

Note that $S$ contains at least $k-1$ vertices. 
By hypothesis, for a fixed $S$, at most $r$ of those $d_s^i$ are non-zero;
say $d_S^{i_1}, \ldots, d_S^{i_r}$. Applying the convex inequality
\eqref{eq:convexity}, we have
$$\sum_{i\in I} {d_S^i\choose s_{l+1}} =
\sum_{j=1}^r {d_S^{i_j}\choose s_{l+1}}\geq r {d_S/r \choose s_{l+1}},$$
provided $d_S >r(s_{l+1}-1)$. 

Let $\bar d_l$ be the average of $d_S$. 
By equation \eqref{eq:si} and inductive hypothesis, we have
\begin{equation}
  \label{eq:dl}
\bar d_l\geq \frac{\sum_{S}\sum_{i\in I} d_S^i}{m^{k -l-1+ \sum_{j=1}^ls_j}}
\geq  \left(1+o(1)\right)
 \frac{m\left(\frac{C}{m^\delta}\right)^{\prod_{j=1}^l s_j}} 
{\prod_{j=1}^l \left(s_j! r^{s_j-1}\right)^{\prod_{u=j+1}^ls_u}}.
\end{equation}

Let $\Sa$ be the set of $S$ satisfying $d_S >r(s_{l+1}-1)$.
Let $\bar d^*$ be the average of $d_S$ over $S\in \Sa$.
Clearly,  $\bar d_l^*\geq \bar d_l$ since $\bar d_l\gg r(s_{l+1}-1)$
Thus,
\begin{align*}
f(s_1,\ldots, s_l, s_{l+1},1,\ldots,1)
&\geq\sum_{S\in \Sa}\sum_{i\in I} {d_S^i\choose s_{l+1}}\\
&\geq  \sum_{S\in \Sa}  r {d_S/r \choose s_{l+1}}\\
&\geq r|\Sa|{\bar d^*_l/r\choose s_{l+1}}\\
&=\frac{|\Sa| \bar d^*_l}{s_{l+1}} {d^*_l/r\choose s_{l+1}-1}\\
&\geq \frac{ (\bar d_l -r(s_{l+1}-1))m^{k-1 +\sum_{j=1}^l(s_j-1)}}{s_{l+1}} 
{\bar d_l/r\choose s_{l+1}-1}\\
&= \left(1+O\left(\frac{1}{\bar d_l}\right)\right)\frac{\bar d_l^{s_{l+1}}}{s_{l+1}! r^{s_{l+1}-1}} m^{k-1 +\sum_{j=1}^l(s_j-1)}.
\end{align*}
Combining with equation \eqref{eq:dl}, we get
$$f(s_1,\ldots, s_l, s_{l+1},1,\ldots,1)
\geq \left(1+o(1)\right)
\frac{\left(\frac{C}{m^\delta}\right)^{\prod_{j=1}^{l+1} s_j}} 
{\prod_{j=1}^{l+1} \left(s_j! r^{s_j-1}\right)^{\prod_{u=j+1}^{l+1}s_u}}
m^{k -l-1+ \sum_{j=1}^{l+1}s_j}.
$$
The inductive proof is finished.

Applying Claim (a) with $l=k-1$, we get
\begin{align} \nonumber
f(s_1,s_2,\ldots, s_{k-1},1) &\geq
\left(1+o(1)\right)
\frac{\left(\frac{C}{m^\delta}\right)^{\prod_{j=1}^{k-1} s_j}} 
{\prod_{j=1}^{k-1} \left(s_j! r^{s_j-1}\right)^{\prod_{u=j+1}^{k-1}s_u}}
m^{1+ \sum_{j=1}^{k-1}s_j} \\
&=\left(1+o(1)\right) \frac{C^{\prod_{j=1}^{k-1} s_j} m^{\sum_{j=1}^{k-1}s_j}}
{\prod_{j=1}^{k-1} \left(s_j! r^{s_j-1}\right)^{\prod_{u=j+1}^{k-1}s_u}}.
\label{eq:lb}
\end{align}
For any $S\in {V_1\choose s_1}\times \cdots \times {V_{k-1}\choose s_{k-1}}$, 
let $d_S$ be the 
number of vertices $v$ in $V_{l+1}$ such that the edges in the induced subgraph
of $H'$ on $S\times \{v\}$ are monochromatic. Since $H'$ contains no
monochromatic copy of $K^{(k)}_k(s_1,\ldots, s_k)$, we have
$d_S\leq rs_k$. It implies
\begin{equation}
  \label{eq:ub}
f(s_1,s_2,\ldots, s_{k-1},1)=\sum_S d_S\leq rs_k m^{\sum_{j=1}^{k-1}s_j}.  
\end{equation}
Choosing $C>2(rs_k)^{\frac{1}{\prod_{u=1}^{k}s_u}}\cdot
\prod_{j=1}^{k-1} \left(s_j! r^{s_j-1}\right)^{\frac{1}{\prod_{u=1}^{j}s_u}}$,
equations \eqref{eq:lb} and \eqref{eq:ub} contradict  each other.
Hence, $|E(H')|<C m^{k-\delta}$. It implies $\sum_{i\in I}|E(H_i)|=O(m^{k-\delta})
=O(n^{k-\delta})$.
The proof of the lemma is finished.
\hfill $\square$

\section{Proof of main Theorem}
We  need the following two lemmas  on  binomial  coefficients. 
\begin{lemma}\label{3.1}
 (see Lemma 2.1 of \cite{GriLu})
For any positive integer $n$, we have
  \begin{equation}
\label{eq:11}
 \frac{1}{2^n} \sum_{|i-\frac{n}{2}|> 2\sqrt{n\ln n}} {n\choose i}< \frac{2}{n^2}.
  \end{equation}
\end{lemma}

\begin{lemma} \label{binom} 
For any $i,j\in (\frac{n}{2}-2\sqrt{n\ln
n},\frac{n}{2}+2\sqrt{n\ln n})$, if $|i-j|=o \left(\frac{\sqrt{n}}{\sqrt{\ln n}}\right)$,
then 
\begin{equation}
  \label{eq:ij}
\frac{{n\choose i}}{{n\choose j}}=1+o(1).  
\end{equation}
 \end{lemma}
{\bf Proof:} Without loss of generality, we can assume $j>i\geq \frac{n}{2}$.
 We have
\begin{align*}
 \frac{{n\choose i}}{{n\choose j}} &= \prod_{l=1}^{j-i}  \frac{{n\choose i+l-1}}{{n\choose i+l}}\\
&= \prod_{l=1}^{j-i} \frac{i+l}{n-i-l+1}\\
&= \prod_{l=1}^{j-i}  \left( 1+ \frac{2(i+l)-n-1}{n-i-l+1} \right).
\end{align*}
Since $i+l\in (\frac{n}{2}-2\sqrt{n\ln
n},\frac{n}{2}+2\sqrt{n\ln n})$, we have

$$ \frac{|2(i+l)-n-1|}{n-i-l+1} \leq \frac{4\sqrt{n\ln n}+1}{\frac{n}{2}-2\sqrt{n\ln n}}
= (1+o(1)) \frac{8\sqrt{\ln n}}{\sqrt{n}}.$$
Thus, we get
$$ \frac{{n\choose i}}{{n\choose j}} \leq 
\left(1+  (1+o(1)) \frac{8\sqrt{\ln n}}{\sqrt{n}}\right)^{j-i}
=1+o(1).$$
The proof of the lemma is finished.
\hfill $\square$

To prove Theorem \ref{main}, we need to show that for any $\epsilon>0$
any family $\F\subset 2^{[n]}$ of size $(1+\epsilon)\nchn$ must
contain the subposet $P$. Without loss of generality, we can assume
that $\F$ only contains subsets of sizes in the interval
$(\frac{n}{2}-2\sqrt{n\ln n},\frac{n}{2}+2\sqrt{n\ln n})$. This is
because the number of subsets of size not in $I$ 
(see Lemma \ref{3.1}) is at most
$$\sum_{|l-\frac{n}{2}|>2\sqrt{n\ln n}}{n\choose l}\leq
\frac{2^{n+1}}{n^2}=O\left(\frac{\nchn}{n^{3/2}}\right),$$ which is
negligible compared to $\epsilon \nchn$.

Taking a random permutation $\sigma$ of the set $[n]$,
a (random) full chain is the chain
$$\emptyset \subset \{\sigma(1)\}
\subset  \{\sigma(1), \sigma(2)\}
\subset \cdots \subset [n].$$
Let $X$ be the number
of subsets in both $\F$ and a random full chain. 
The expected value of $X$ is exactly the Lubell value of $\F$:
\begin{equation}
  \label{eq:lubell}
\E(X)=h_n(\F)=\sum_{F\in \F}\frac{1}{{n\choose |F|}}.
\end{equation}

It is clear that
\begin{equation}
  \label{eq:ex}
\E(X)\geq \frac{|\F|}{\nchn}=1+\epsilon.  
\end{equation}
 Combining equation \eqref{eq:ex} and 
 the convexity inequality \eqref{eq:convexity} with $s=2$, we have
 \begin{equation}
   \label{eq:ex2lb}
 \E{X\choose 2} \geq  {\E(X)\choose 2}\geq
\frac{\epsilon}{2}\E(X).  
 \end{equation}

%To get the contradiction,  we will show $\E{X\choose 2}=o(\E(X))$.
For any two subsets $A\subseteq B$,
the probability that a random full chain hits both $A$ and $B$ is
$\frac{|A|! (|B|-|A|)!(n-|B|)!}{n!}$. By linearity, we get
\begin{equation}
  \label{eq:ex2}
\E{X\choose 2}=\sum_{\stackrel{A,B\in \F}{A\subset B}}
\frac{|A|! (|B|-|A|)!(n-|B|)!}{n!}.
\end{equation}

The following Lemma was implicitly proved when Griggs
and Lu \cite{GriLu} proved $\La(T)=(1+o(1))\nchn$ for any tree poset of height 2.
The statement works for any poset of height $2$, not just those
having $k$-partite representation. We state it here as a lemma  
for the future references, and also provide a proof for completeness.
\begin{lemma}\label{AB}
Let $P$ be a finite poset of height $2$ 
and $\F$ be a $P$-free $\F$ family of subsets of $[n]$ with the Lubell value
$h_n(\F)\geq 1+\epsilon$. Suppose that every subset in $\F$ has size in the interval
$(\frac{n}{2}-2\sqrt{n\ln n},\frac{n}{2}+2\sqrt{n\ln n})$.
Then, we have
\begin{equation}
  \label{eq:AB}
\sum_{\stackrel{A,B\in \F, |B|-|A|=1}{A\subset B}}
\frac{|A|! (|B|-|A|)!(n-|B|)!}{n!}
\geq (1+o(1))\epsilon h_n(\F).
\end{equation}
\end{lemma}
{\bf Proof:}
Let $Y$ be the random variable
counting a triple $(A,S,B)$ (on the random full chain)
satisfying
\[A\subset S\subset B\quad A,B\in \F.\]
We have
\begin{eqnarray}
\nonumber
  \E(Y)&=& \sum_{\stackrel{A,B\in \F, S}{A\subset S\subset B}}
\frac{|A|! (|S|-|A|)!(|B|-|S|)!(n-|B|)!}{n!} \\
\nonumber
&=& \sum_{\stackrel{A,B\in \F}{A\subset B}}
\frac{|A|! (|B|-|A|)!(n-|B|)!}{n!} \sum_{S\colon A\subset S\subset B}
 \frac{1}{{|B|-|A|\choose |S|-|A|} }\\
\nonumber
&=& \sum_{\stackrel{A,B\in \F}{A\subset B}}
\frac{|A|! (|B|-|A|)!(n-|B|)!}{n!} (|B|-|A|-1) \\
\label{eq:41}
&\geq&  \sum_{\stackrel{A,B\in \F}{A\subset B, |B|-|A|>1}}
\frac{|A|! (|B|-|A|)!(n-|B|)!}{n!}.
\end{eqnarray}

Any poset $P$ of height $2$ is a subposet of $K_{r,r}$ (the complete
height-$2$-poset) for some $r$.
Since $\F$ is $P$-free, there are no
$2r$ subsets $A_1,A_2,\ldots, A_r, B_1,\ldots, B_r\in \F$
satisfying $A_i\subset S\subset B_j$ for $1\leq i\leq r$ and
$1\leq j\leq r$.

For any fixed subset $S$, either ``at most $r-1$ subsets in $\F$
are supersets of $S$'' or ``at most $r-1$ subsets in $\F$ are
subsets of $S$''. Define
\[\G_1=\{ S\mid |S|\in
(\frac{n}{2}-2\sqrt{n\ln n},\frac{n}{2}+2\sqrt{n\ln n}),
\mbox{$S$  has at most
$r-1$ subsets in $\F$}\}.\]
\[\G_2=\{ S\mid |S|\in
(\frac{n}{2}-2\sqrt{n\ln n},\frac{n}{2}+2\sqrt{n\ln n}),
\mbox{$S$  has at most
$r-1$ supersets in $\F$}\}.\]
The union $\G_1\cup \G_2$ covers all subsets with sizes in
$(\frac{n}{2}-2\sqrt{n\ln n},\frac{n}{2}+2\sqrt{n\ln n})$.
Rewrite $\E(Y)$ as
\begin{equation}
  \label{eq:42}
\E(Y) =
\sum_{S\colon ||S|-\frac{n}{2}|<2\sqrt{n\ln n}}
\frac{1}{{n\choose |S|}}
 \sum_{\stackrel{A\in \F}{A\subset S}}
\frac{1}{{|S|\choose |A|}}
 \sum_{\stackrel{B\in \F}{S\subset B}}
\frac{1}{{n-|S|\choose n-|B|}}.
\end{equation}
For $S\in\G_1$, we have
\begin{equation}
  \label{eq:43}
\sum_{B\in \F, S\subset B}\frac{1}{{n-|S|\choose n-|B|}}
\leq \frac{r-1}{\frac{n}{2}-2\sqrt{n\ln n}}=O(\frac{1}{n}).
\end{equation}
It implies
$$  \sum_{S\in \G_1}
\frac{1}{{n\choose |S|}}
 \sum_{\stackrel{A\in \F}{A\subset S}}
\frac{1}{{|S|\choose |A|}}
 \sum_{\stackrel{B\in \F}{S\subset B}}
\frac{1}{{n-|S|\choose n-|B|}}
\leq
 \sum_{S\in \G_1}
\frac{1}{{n\choose |S|}}
 \sum_{\stackrel{A\in \F}{A\subset S}}
\frac{1}{{|S|\choose |A|}}
O\left(\frac{1}{n}\right)
.$$
Recall $\E(X)=h_n(\F)= \sum_{\stackrel{A\in \F}{A\subset S}}
\frac{1}{{|S|\choose |A|}}$ and
$\sum_{S\in \G_1}
\frac{1}{{n\choose |S|}} \leq 4\sqrt{n\ln n}$.
We have
$$
 \sum_{S\in \G_1}
\frac{1}{{n\choose |S|}}
 \sum_{\stackrel{A\in \F}{A\subset S}}
\frac{1}{{|S|\choose |A|}}
 \sum_{\stackrel{B\in \F}{S\subset B}}
\frac{1}{{n-|S|\choose n-|B|}}
\leq O\left( \frac{\sqrt{\ln n}}{\sqrt{n}}\E(X)\right).
$$

Similarly, we have
$$
  \sum_{S\in \G_2}
\frac{1}{{n\choose |S|}}
 \sum_{\stackrel{A\in \F}{A\subset S}}
\frac{1}{{|S|\choose |A|}}
 \sum_{\stackrel{B\in \F}{S\subset B}}
\frac{1}{{n-|S|\choose n-|B|}}
=O\left( \frac{\sqrt{\ln n}}{\sqrt{n}}\E(X)\right).
$$
Thus, we have
\begin{equation}
\label{eq:y}
\E(Y)=O\left( \frac{\sqrt{\ln n}}{\sqrt{n}}\E(X)\right)
=o(\epsilon \E(X)).
\end{equation}
 
Combining 
inequalities \eqref{eq:ex2lb}, \eqref{eq:41}, \eqref{eq:y}, 
with equation \eqref{eq:ex2},  we have
\begin{equation}
  \label{eq:46}
\sum_{\stackrel{A,B\in \F, |B|-|A|=1}{A\subset B}}
\frac{|A|! (|B|-|A|)!(n-|B|)!}{n!}
= \E{X\choose 2}-\E(Y)
\geq (1-o(1))\epsilon\E(X).
\end{equation}
The proof of Lemma is finished.
\hfill$\square$

% In the rest of proof, we require that $P$ has a $k$-partite represenation.
% we will estimate the number of $k$-configurations, which is
% essentially the $k$-dimensional sublattice with the minimum element
% $S$ and the maximum element $B$ such that all elements of the top 2
% levels belong in $\F$. By the average argument, there are many
% $k$-configurations with the same minimum element $S$. Then we
% construct a collection of $k$-uniform hypergraphs, which satisfies
% the conditions of Lemma \ref{cturan}. Each $k$-uniform edge in this
% family is corresponding to a $k$-configuration. By Lemma \ref{cturan},
% the number of $k$-configurations is small; and we get the desired
% contradiction. Here is the detail proof.

{\bf Proof of Theorem \ref{main}:} Now we assume that $P$ has a
$k$-partite representation and $\F$ is a $P$-free $\F$ family of
subsets of $[n]$ with the Lubell value $h_n(\F)=
1+\epsilon$. We further assume that every subset in $\F$ has size in the
interval $(\frac{n}{2}-2\sqrt{n\ln n},\frac{n}{2}+2\sqrt{n\ln n})$.
Let $X$ be the random variable couting the number of subsets of $\F$
hit by a random full chain. Note $\E(X)=h_n(\F)$.
By Lemma \ref{AB}, we have
\begin{equation}
    \label{eq:47}
\sum_{\stackrel{A,B\in \F, |B|-|A|=1}{A\subset B}}
\frac{|A|! (|B|-|A|)!(n-|B|)!}{n!}
\geq (1-o(1))\epsilon\E(X).
\end{equation}
We define $N(B)=\{A\in \F\mid  A\subset B, |A|=|B|-1\}$
and $d(B)=|N(B)|$.
We have
\begin{equation}
  \label{eq:db}
  \sum_{B\in \F}\frac{1}{{n\choose B}}\frac{d(B)}{|B|}
=\sum_{\stackrel{A,B\in \F}{A\subset B,|B|-|A|=1}}
\frac{|A|! (|B|-|A|)!(n-|B|)!}{n!}.
\end{equation}
%Equation \eqref{eq:db} and inequality \eqref{eq:46} imply
%that there are many $B\in \F$ with $d(B)\geq \frac{\epsilon}{2}n$.
%Let $\B:=\{B\in \F\colon d(B)>\frac{\epsilon}{2}n\}$.

Let $\bar d:= \frac{1}{\E(X)}\sum_{B\in \F}\frac{d(B)}{{n\choose B}}$
be the weighted average of $d(B)$.
Since $|B|=(1+o(1))\frac{n}{2}$ for any $B\in \F$,
by equation \eqref{eq:db} and inequality \eqref{eq:47}, 
we have
\begin{align*}
 \bar d &= \frac{1}{\E(X)}\sum_{B\in \F}\frac{d(B)}{{n\choose B}} \\
&= (1+o(1))\frac{n}{2\E(X)} \sum_{B\in \F}\frac{d(B)}{{n\choose B}|B|} \\
&= (1+o(1))\frac{n}{2\E(X)} \sum_{\stackrel{A,B\in \F}{A\subset B,|B|-|A|=1}}
\frac{|A|! (|B|-|A|)!(n-|B|)!}{n!}\\
&\geq  (1+o(1))\frac{\epsilon n}{2}.
\end{align*}

A pair of sets $(S, B)$ is said to form a {\em $k$-configuration} if 
\begin{enumerate}
 \item $S\subset B$, $|S|=|B|-k$, and $B\in \F$;
\item for any $x\in B\setminus S$,
$B\setminus \{x\} \in \F$ .
\end{enumerate}

Since $||B|-\frac{n}{2}|\leq 2\sqrt{n\ln n}$,
 $|S|$ belongs to the interval 
$J:=(\frac{n}{2}-2\sqrt{n\ln n}-k,\frac{n}{2}+2\sqrt{n\ln n}-k)$.
Set ${\cal S}:=\cup_{s\in J}{[n]\choose s}$.
For any $S\in \cal S$,
let $L(S)$ be the number of such configurations  over a fixed set $S$.
We have
\begin{align*}
\sum_{S\in \cal S}\frac{L(S)}{{n\choose |S|}} 
&= (1+o(1))\sum_{S\in \cal S}  \frac{L(S)}{{n\choose |S|+k}}\\
&=  (1+o(1))\sum_{B\in \F} \frac{1}{{n\choose |B|}}{d(B)\choose k}\\
&\geq (1+o(1))\E(X) {\bar d \choose k}\hspace*{2cm} \mbox{ by the convexity inequality \eqref{eq:convexity}}
\\
&\geq (1+o(1))\frac{\epsilon^k}{2^k k!} n^k \E(X)\\
&\geq \frac{\epsilon^k}{2^k k!} n^k.
\hspace*{2cm} 
\end{align*}
Partition $J$ into small sub-intervals $\{J_\lambda\}_{\lambda\in \Lambda}$ with equal length  
$\frac{\sqrt{n}}{\ln n}$.
There are  $4\ln^{3/2} n$ of such sub-intervals. Setting
${\cal S}_\lambda:=\cup_{s\in J_\lambda}{[n]\choose s}$, we have
${\cal S}=\cup_{\lambda\in \Lambda}{\cal S}_\lambda$. By an average argument,
there is a $\lambda_0\in \Lambda$ so that
\begin{equation}
  \label{eq:lambda0}
\sum_{S\in {\cal S}_{\lambda_0}}\frac{L(S)}{{n\choose |S|}}
\geq \frac{1}{4\ln^{3/2} n} \sum_{S\in \cal S}\frac{L(S)}{{n\choose |S|}} 
\geq \frac{\epsilon^k}{2^{k+2} k!}\frac{n^k}{\ln^{3/2} n}.  
\end{equation}

Suppose that ${n\choose s}$ for $s\in J_{\lambda_0}$ reaches the maximum at $s=s_0$.
Note that $|s-s_0|\leq \frac{\sqrt{n}}{\ln n}$.
By Lemma \ref{binom}, we have
\begin{equation}
  \label{eq:s0}
{n\choose s}=(1-o(1)){n\choose s_0}.  
\end{equation}

Combining equations \eqref{eq:lambda0} and \eqref{eq:s0},
we get
\begin{equation}
  \label{eq:21}
 \sum_{S\in {\cal S}_{\lambda_0}}L(S)\geq (1-o(1)) 
\frac{\epsilon^k}{2^{k+2} k!}\frac{n^k}{\ln^{3/2} n}{n\choose s_0}.
\end{equation}
Observe that there is a chain decomposition of ${\cal S}_{\lambda_0}$ into ${n\choose s_0}$
chains. There exists one chain $\C$ satisfying
\begin{equation}
  \label{eq:22}
  \sum_{S\in \C}L(S)\geq \frac{1}{{n\choose s_0}} \sum_{S\in {\cal S}_{\lambda_0}}L(S)
\geq (1-o(1)) 
\frac{\epsilon^k}{2^{k+2} k!}\frac{n^k}{\ln^{3/2} n}.
\end{equation}
For any $S\in \C$, we define a $k$-uniform hypergraph $H_S$ on the vertex set $[n]$
as follows: a $k$-set $F$ is an edge of $H_S$ if $S\cap F=\emptyset$ and
$(S, S\cup F)$ forms a $k$-configuration.
% \begin{enumerate}
% \item $S\cap \{i_1,i_2,\ldots, i_k\}=\emptyset$.
% \item $S\cup \{i_1,i_2,\ldots, i_k\}\in \F$.
% \item For $1\leq l\leq k$, $S\cup \{i_1,i_2,\ldots, i_k\}\setminus\{i_l\}\in \F$.
% \end{enumerate}

Let $\Pa\subset {[l]\choose k-1} \cup {[l]\choose k}$ be 
the $k$-representation of $P$ and  $G(\Pa)$ be the
$k$-uniform hypergraph associated with $\Pa$. 
Since $G(\Pa)$ is $k$-partite, there is a $k$-partition  
$$[l]=V_1\cup  V_2 \cdots \cup V_k$$
such that all edges of $G(\Pa)$ are crossing.
For $1\leq i\leq k$, set $s_i:=|V_i|$.
Clearly, we have   $G(\Pa) \subset K^{(k)}_k (s_1, \ldots, s_k)$.

{\bf Claim b:} The hypergraph $H_S$ contains no copies of
 $K^{(k)}_k (s_1, \ldots, s_k)$ as a sub-hypergraph.
Otherwise, $H_S$ contains $G(\Pa)$ as a subgraph. 
By the definition of $H_S$, $\{S\cup F\}_{F\in \Pa}\subset \F$.
As a poset,  $\{S\cup F\}_{F\in \Pa}\subset \F$ is isomorphic to $\Pa$.
Thus,  $\F$ contains a subposet $P$.

{\bf Claim c:} For any $(k-1)$-set $T$, the number of edges 
of $H_S$ (for $S\in \C$) containing $T$ is at most $r$. Otherwise, 
there exists a chain
$$S_1\subset S_2 \subset \cdots \subset S_r$$
such that $T\in E(H_{S_i})$ for all $1\leq i \leq r$.
By the definition of $H_{S_i}$, we have $T\cup S_i\in \F$.
Thus, $(S_1\cup T),(S_2\cup T), \cdots, (S_r\cup T)$ forms
an $r$-chain in $\F$. This chain contains the subposet $P$.
Contradiction.

By Claims (b) and (c), the collection $\Ha:=\{H_S\}_{S\in \C}$ satisfies
the conditions of Lemma \ref{cturan}. Hence,
the total number of edges in $\Ha$ is $O(n^{k-\delta})$,
where $\delta =\left(\prod_{i=1}^{k-1}s_i\right)^{-1}$ is a positive constant.
Note that an edge in $H_S$ is 1-1 corresponding to a $k$-configuration
$(S,B)$. Thus, we have
\begin{equation}
  \label{eq:hs}
   \sum_{S\in \C}L(S)=O(n^{k-\delta}).
\end{equation}

Combining equations \eqref{eq:22} with \eqref{eq:hs},
we get $$\epsilon^k=O\left(\frac{\ln^{3/2} n}{n^\delta}\right)=o(1).$$
This contradicts the assumption that $\epsilon$ is a constant.
The proof of the theorem is finished. \hfill $\square$


\begin{thebibliography}{99}
\bibitem{AKT}
N.~Alon, A.~Krech and T.~Szab\'o, Tur\'an’s theorem in the hypercube, 
{\em SIAM J. Discrete Math.} 21 (2007), 66--72.

\bibitem{ARSV} N.~Alon, R.~Radoi\v{c}\'ic, B.~Sudakov and
  J.~Vondr\'ak, A Ramsey-type result for the hypercube, {\em J. Graph
  Theory} {\bf 53} (2006), 196--208.

%\bibitem{AxeMar} M. Axenovich and R. Martin, A note on short cycles in the hypercube, 
%{\em Discrete Math.} {\bf 306} (2006), 2212--2218.
%
\bibitem{AxeManMar} M.~Axenovich, J.~Manske, and R.~Martin,
{$Q_2$-free families in the Boolean lattice}, {\em Order}
published online: 15 March 2011.

\bibitem{Buk}
B.~Bukh, {Set families with a forbidden poset}, {\em Elect. J.
Combin.} {\bf 16} (2009), R142, 11p.

%\bibitem{CarKat}
%T. Carroll and G. O. H. Katona, Bounds on
%  maximal families of sets not containing three sets with
%$A\cup B\subset C, A \not\subset B$, {\em Order} {\bf 25} (2008)
%229--236.
%
\bibitem{chung}
F.~Chung, Subgraphs of a hypercube containing no small even cycles, {\em J. Graph Theory}
{\bf 16} (1992), 273--286.

%\bibitem{conder}
%M. Conder, Hexagon-free subgraphs of hypercubes, {\em J. Graph Theory} 
%{\bf 17} (1993), 477--479.
%
\bibitem{Conlon}
D.~Conlon, An extremal theorem in the hypercube. {\it Electron. J. Combin.} 
{\bf 17} (2010), \#R111.


\bibitem{DebKat} A.~De~Bonis and G.~O.~H.~Katona, {Largest families
without an $r$-fork}, {\em Order} {\bf 24} (2007), 181--191.

\bibitem{DebKatSwa}
A.~De~Bonis, G.~O.~H.~Katona and K.~J.~Swanepoel, {Largest family
without $A\cup B \subset C\cap D$}, {\em J. Combin. Theory (Ser.
A)} {\bf 111} (2005), 331--336.


\bibitem{Erdos45} P.~Erd\H{o}s,
On a lemma of Littlewood and Offord, {\em Bull. Amer. Math. Soc.}
{\bf 51} (1945), 898--902.

\bibitem{Erdos64} P.~Erd\H{o}s, On extremal problems of graphs and generalized graphs, 
{\em Israel J. Math.} {\bf 2} (1964), 183--190.

\bibitem{c14} Z.~F\"uredi and L.~\"Ozkahya, On 14-cycle-free subgraphs of the hypercube, 
{\em Combin. Probab. Comput.} {\bf 18} (2009), 725--729.

\bibitem{c14+} Z.~F\"uredi and L.~\"Ozkahya, On even-cycle-free subgraphs of the hypercube, 
{\em Electronic Notes in Discrete Mathematics} {\em 34} (2009), 515--517.


%\bibitem{GreKle} C. Greene and D. J. Kleitman, {Proof techniques
%in the theory of finite sets}, in: G.C. Rota (ed.), {\em Studies
%in Combinatorics}, MAA Studies in Mathematics {\bf 17}, MAA,
%Providence, 1978,  pp. 22--79.
%
%\bibitem{GriKat} J. R. Griggs and G. O. H. Katona, {No
%four subsets forming an $N$}, {\em J. Combinatorial Theory (Ser.
%A)} {\bf 115} (2008), 677--685.
%
%\bibitem{GriLi} J. R. Griggs and W.-T. Li, {The partition method
%for poset-free families}, preprint (2011).
%
\bibitem{GriLi2} J.~R.~Griggs and W.-T.~Li, {Uniformly L-bounded
posets}, preprint (2011).

\bibitem{GriLiLu} J.~R.~Griggs, W.-T.~Li,  and L.~Lu,
Diamond-free Families, {\em Journal of Combinatorial Theory
Ser. A}, {\bf 119} (2012) 310-322.

\bibitem{GriLu} J.~R.~Griggs and L.~Lu, {On families of subsets
with a forbidden subposet}, {\em Combinatorics, Probability, and
Computing} {\bf 18} (2009), 731--748.

%\bibitem{Kat08} Gyula O. H. Katona, Forbidden inclusion patterns in the
%  families of subsets (introducing a method), in {\em Horizons of
%    Combinatorics}, Bolyai Society Mathematical Studies, {\bf 17}, Bolyai
%  Mathematical Society, Budapest and Springer-Verlag, 2008,
%  pp. 119--140.
%
\bibitem{KatTar}  G.~O.~H.~Katona and T.~G.~Tarj\'an, {Extremal problems
with excluded subgraphs in the $n$-cube},  in: M. Borowiecki, J.
W. Kennedy, and M. M. Sys\l o (eds.) {\bf Graph Theory}, \L
ag\'ow, 1981, {\em Lecture Notes in Math.}, {\bf 1018} 84--93,
Springer, Berlin Heidelberg New York Tokyo, 1983.

\bibitem{KMY} L.~Kramer, R.~Martin, M.~Young,
On diamond-free subposets of the Boolean lattice,
\verb+http://arxiv.org/abs/1205.1501+.


%\bibitem{WTL} W.-T. Li, {Extremal Problems on Families of Subsets with
%Forbidden Subposets}, Ph.D. dissertation, University of South
%Carolina, 2011.
%
%\bibitem{Lub}  D. Lubell, {A short proof of Sperner's lemma},
%{\em J. Combin. Theory} {\bf 1}(1966), 299.
%
\bibitem{Spe} E.~Sperner,
{Ein Satz \"{u}ber Untermegen einer endlichen Menge}, {\em Math.
Z.} {\bf 27} (1928), 544--548.

\bibitem{Sta} R.~P.~Stanley, {\em Enumerative Combinatorics Vol. 1},
Cambridge University Press, 1997.

\bibitem{Tha} H.~T.~Thanh,
An extremal problem with excluded subposets in the Boolean lattice,
{\em Order} {\bf 15} (1998), 51--57.

\end{thebibliography}
\end{document}